\documentclass{amsart}
\usepackage{amsmath, amstext, amsbsy}
\usepackage{amssymb}

\newtheorem{theorem}{Theorem}[section]
\newtheorem{lemma}[theorem]{Lemma}
\newtheorem{proposition}{Proposition}[section]   

\newtheorem{cor}{Corollary}[section]

\theoremstyle{definition}
\newtheorem{definition}[theorem]{Definition}
\newtheorem{example}[theorem]{Example}

\theoremstyle{remark}
\newtheorem{remark}[theorem]{Remark}

\numberwithin{equation}{section}

\newcommand{\la}{\lambda}
\newcommand{\vi}{e^{i\alpha/2}}

\begin{document}


\title[Symmetric polynomials and $U_q(\widehat{sl}_2)$]
    {Symmetric polynomials and $U_q(\widehat{sl}_2)$
    }
\author{Naihuan Jing}
\address{Department of Mathematics,
   North Carolina State University,
   Raleigh, NC 27695-8205,
   U.S.A.}
\email{jing@math.ncsu.edu}
\thanks{Research supported in part by NSA grant
MDA904-97-1-0062 and Mathematical Sciences Research Institute.}


\date{February 17, 1999; Revised version: June 22, 1999}


\begin{abstract}
We study the explicit formula of 
Lusztig's integral forms of the level one 
quantum affine algebra $U_q(\widehat{sl}_2)$ in the 
endomorphism ring of symmetric functions in infinitely many variables tensored
with the group algebra of $\mathbb Z$. 
Schur functions are 
realized as certain orthonormal basis vectors in 
the vertex representation associated to the standard Heisenberg algebra. 
In this picture the Littlewood-Richardson rule is expressed by integral 
formulas, and is used to define the action of 
Lusztig's $\mathbb Z[q, q]$-form of $U_q(\widehat{sl}_2)$ on Schur polynomials. 
\end{abstract}

\maketitle

\section{Introduction}\label{S:intro}

 The relation between vertex representations and the symmetric functions is one
of the interesting aspects of affine Kac-Moody algebras, quantum affine 
algebras
and the vertex operator algebra. In the early 1980's 
the Kyoto school \cite{DJKM}
found that the polynomial solutions of KP hierarchies are obtained by 
Schur polynomials. This breakthrough was achieved in their work
of describing the KP and KdV hierarchies in terms of affine Lie algebras.
On the other hand I. Frenkel \cite{F1}
identified the
two constructions of the affine Lie algebras via the vertex operators, which 
put the boson-fermion correspondence
in a rigorous formulation. Later I.~Frenkel further
used this idea \cite{F2} to find out that the boson-fermion correspondence can
also be used to realize the Schur polynomials in the study of 
irreducible characters of the symmetric
group $S_n$. Schur functions also appeared in Lepowsky and Primc
construction \cite{LP} of certain bases for higher level representations of
the affine Lie algebra $\widehat{sl}(2)$. 

In \cite{J1, J2} the author advanced the vertex operator approach to the 
symmetric functions and realized 
Schur's Q-functions, which assume the same role of Schur functions
in the context of (non-trivial) irreducible characters of the double 
covering groups of $S_n$, and more generally  
Hall-Littlewood polynomials were treated in the same way. 
These families of symmetric polynomials appear naturally as
orthogonal bases in the vertex representation. 
In \cite{J1} I. Frenkel's original work
on $S_n$ was reviewed together with Bernstein's formula (see \cite{Z}).
The exact formulation of Hall-Littlewood polynomials in terms of
the boson-fermion correspondence appeared later in the author's work
\cite{J5}.
Since Macdonald's seminal extension \cite{M}
of Hall-Littlewood polynomials
 it immediately
brought up the question of the vertex operator approach to
the more general Macdonald polynomials. 

In \cite{J4} we also found that the vertex operator realization of certain
Macdonald polynomials is governed by the basic hypergeometric
functions of type $_4\phi_3$.
It was realized \cite{J3} that
the complete vertex operator realization of Macdonald polynomials
is related to several new questions, among them the most interesting
question was how one could realize the Macdonald polynomials in the vertex
representations \cite{FJ} of the quantum affine algebra $U_q(\widehat{sl}_2)$.
In this case the 
underlying symmetry of the analogous family of
symmetric functions are expressed by
certain  infinite dimensional quadratic algebras generalizing
the Clifford-Weyl algebra.

Recently J. Beck, I. Frenkel and the author \cite{BFJ}
have used vertex operators to study the canonical bases 
for the level one irreducible modules for the quantum affine algebra 
$U_q(\widehat{sl}_2 )$.
Macdonald polynomials (zonal case) are seen as some "canonical" bases sitting
between Kashiwara and Lusztig's canonical and dual canonical bases 
\cite{L, K}. This
essentially answered the 
question about the vertex realization of Macdonald polynomials.
The Macdonald basis constructed in \cite{BFJ} also
satisfy the characteristic properties of bar invariance and orthogonality 
under the
Kashiwara form. 
The transition matrix from the canonical basis to Macdonald 
basis is triangular,
integral and bar-invariant. We also conjectured its positivity.
Since the transition matrix from Macdonald
polynomials to (modified) 
Schur polynomials is also triangular, it is natural to ask how
does the quantum affine algebra act on Schur polynomials. 
 
Our first goal in this paper is to answer the question of realizing
the quantum affine algebra $U_q(\widehat{sl}_2)$ by Schur functions. We give
explicit formulas to realize the level one representation with the help of
Littlewood-Richardson rule.
 In order to do this we first redevelop the vertex approach 
to Schur functions. 
Mixed products of vertex operators and dual vertex operators are expressed 
in terms of Schur functions. 
We showed
that the underlying symmetry of the Young tableaus of the 
various classical symmetric functions are governed by
certain Clifford-Weyl type algebras, where the simplest case
of $S_n$-symmetry or the linkage symmetry
corresponds to the infinite dimensional Clifford algebra.

Using  this 
the quantum affine algebra $U_q(\widehat{sl}_2 )$ 
is represented as $\mathbb Z[q,q^{-1}]$-linear operators on
the ring of symmetric functions. Moreover, 
the lattice of the irreducible module $V(\Lambda_i)$ ($i=0, 1$):
$$
\mathcal L(\Lambda_i)=\Lambda(x_1,\cdots,x_n
\cdots )\bigotimes {\mathbb C}[{\mathbb Z\alpha}]e^{i\alpha/2}
$$
is shown to be an integral lattice for Lusztig's $\mathbb
Z[q, q^{-1}]$-form of 
$U_q(\widehat{sl}_2 )$ (generated by divided
powers). On the other hand Beck, Chari and Pressley
\cite{BCP} have shown that a related lattice 
$\mathcal L(\Lambda_i)'=\Lambda(x_1',\cdots,x_n'
\cdots )\bigotimes {\mathbb C}[{\mathbb Z\alpha}]e^{i\alpha/2}$
is contained in $V(\Lambda_i)_{\mathbb Z[q, q^{-1}]}$,
where $\sum_i {x'}^n=(1-q^{4n})\sum x_i^n$. 
The explicit realization of the Lusztig's integral lattice
$V(\Lambda_i)_{\mathbb Z[q, q^{-1}]}$ will be treated elsewhere.

We further write down the action of
divided powers of current operators
in terms of Littlewood-Richardson rule. Our method stems from the
trick \cite{J3} of expressing Schur functions in terms of
a deformed Heisenberg algebra inside $U_q(\widehat{sl}_2)$. Thus we are able to write down everything explicitly
using the vertex representation, which partly generalizes Garland's
work \cite{G} (see also \cite{CP}).
The method in this paper can also be generalized to quantum affine algebras 
of ADE types, and this may provide more information about Schur 
functions as crystal bases \cite{BCP}. Similarly our formulas might help to understand the positivity 
conjecture \cite{BFJ}.

Our second goal is to give 
a vertex operator approach to the Littlewood-Richard\-son rule. 
In particular an
integral formula for the Littlewood-Richard\-son rule is found. In the
fermionic picture of the crystal basis of quantum affine algebras
Misra and Miwa \cite{MM} used some insertion and deletion procedure on
Young tableaus
to describe the action of quantum affine algebra at $q=0$. The
explicit rule in the fermionic construction \cite{H} suggests
that there should be a corresponding rule in the homogeneous
construction. This  turns out to be an explicit
formula in terms of the Littlewood-Richardson rule in our work (see
Proposition \ref{P:littlewood} and Theorem \ref{theorem3}).
  
Our explicit formulas of the divided powers of Chevalley
generators suggest that there should be corresponding
formulas in the fermionic picture. 
The answer to this question will be helpful to understand the real meaning of Littlewood-Richardson rule in the
boson-fermion correspondence.

The paper is organized as follows. In section \ref{S:schur} we
redevelop our vertex operator approach to Schur polynomials
and express all mixed products of dual vertex operators in terms of Schur functions ($q=1$ case). In section \ref{S:quantum} we first recall the Frenkel-Jing vertex representation of $U_q(\widehat{sl}_2)$
and obtain explicit formulas for the Drinfeld generators
in terms of the Schur basis constructed in Section \ref{S:schur}. 
We then enter the Littlewood-Richardson rule to give the
formulas for the divided powers of the current operators
$X_n^{\pm}$. In the last section (Sec. \ref{S:comb}) we recast the action
in terms of the Chevelley generators, which provides
a simple combinatorial model for the homogeneous
picture of the basic module for $U_q(\widehat{sl}_2)$.
We show that the lattice of Schur bases contains
Lusztig's integral lattice of divided powers. 

\bigskip

{\bf Acknowledgments.} The author wishes to thank I. Frenkel  for helpful 
discussions.  
He also thanks M.L.~Ge
for the warm
hospitality in the summers of 1994 and 1997 at Nankai Institute of 
Mathematics,
where the author had the privilege to stay in the S.S.~Chern Villa
(Home of Geometer) to work out the main computations. 
The author is also grateful to the Mathematical Sciences
Research Institute in Berkeley for providing a lovely
environment in the final stage of this work.

\section{Schur functions and vertex operators}\label{S:schur}

Let $\Lambda_F$ be the ring of symmetric functions in infinitely many 
 variables $x_1, x_2, \cdots $ over the field  $F$. In this section we take $F=\mathbb Q$, and later we will take $F=\mathbb Q(q^{1/2}, q^{-1/2})$. 

A partition $\la=(\la_1, \la_2, \cdots, \la_l)$ of $n$, denoted 
$\lambda\vdash n$,
is a special decomposition of $n$:
$n=\la_1+\cdots+\la_l$ with $\la_1\geq \dots \geq\la_l\geq 1$. $l$ is called the length of $\la$. We will
 identity $(\la_1, \cdots, \la_l)$ with
$(\la_1, \cdots, \la_l, 0, \cdots, 0)$ if we want
to view $\la$ in $\mathbb Z^n$ when
$n\geq l(\la)$. Sometime we
prefer to use another notation for $\la$: $(1^{m_1}2^{m_2}\cdots)$ where $m_i$ is the number of times that $i$ appears among the parts of $\la$. The set of partitions will be denoted by $\mathcal P$.

There are several well-
 known bases in $\Lambda_F$ parameterized by partitions: the power sum symmetric
 functions $p_{\lambda}=p_{\lambda_1}\cdots p_{\lambda_l}$ with
 $p_n=\sum x^n_i$ ($\mathbb Q$-basis), 
the monomial symmetric functions $m_{\lambda}(x_1, \cdots, x_n)=
 \sum_{\sigma}x_1^{\sigma(\lambda_1)}\cdots x_n^{\sigma(\lambda_n)}$, where $\sigma$ runs through distinct permutations of $\lambda$ as tuples; 
and the Schur functions
 $s_{\lambda}$ which are over $\mathbb Z$. In terms of finitely many 
variables Schur function is given by the Weyl character formula:
 \begin{eqnarray}\label{E:weyl}
 s_{\lambda}(x_1, \cdots, x_n)=\frac{\sum_{\sigma\in \mathfrak S_n}sgn(\sigma)x^{
 \sigma(\lambda +\delta)}}{\Pi_{i<j}(x_i-x_j)},
 \end{eqnarray}
where $\delta =(n-1, n-2, \cdots , 1, 0),\,  \lambda =(\lambda_1, \cdots, \lambda_n)$. Here some $\la_i$ may be zero.

Note that this formula also works for any tuple $\mu$. In general
$s_{\mu}=0$ or $(-1)^{l(\sigma)} s_{\lambda}$, where $\lambda=
\sigma(\mu+\delta)-\delta$ for some permutation $\sigma$
and $l(\sigma)$ is the length of the permutation $\sigma$. 
This important property is still true for the Schur function
in infinitely many variables, though there is a less satisfactory 
formula in algebraic
combinatorics in that case. 
We will see that this symmetry property 
is manifested in our vertex operator approach.

 Let's recall the vertex operator approach to Schur functions \cite{J1}.
Let $\{b_n|n\neq 0 \} \bigcup \{c\}$ be the set of generators of the 
Heisenberg algebra with defining relations
\begin{equation}
 [b_m, b_n]=m {\delta}_{m,-n} c, \qquad [c, b_m]=0.
\end{equation}
The Heisenberg algebra has a canonical natural
 representation in the $\mathbb Q$-space
$V=Sym({b_{-n}}'s ), $ the symmetric algebra generated by the $b_{-n}$, 
$n \in {\mathbb N}$. The action is given by
\begin{align}
b_{-n}.v&=b_{-n}v, \qquad
b_n.v=n\frac{\partial v}{\partial b_{-n}}, \\
c.v&=v\, .
\end{align}
It is clear that $1$ is the highest weight vector in $V$.

Let us introduce two vertex operators (cf. $t=0$ in \cite{J2}):

\begin{align}
S(z) &= exp(\sum_{n=1}^{\infty} \frac{b_{-n}}{ n}z^n)
exp(-\sum_{n=1}^{\infty} \frac{b_{ n}}{ n}z^{-n})\\
&=\sum_{n\in\mathbb Z}S_nz^{-n},\nonumber\\
S^*(z) &= exp(-\sum_{n=1}^{\infty} \frac{b_{-n}}{ n}z^n)
exp(\sum_{n=1}^{\infty} \frac{b_{ n}}{ n}z^{ n})\\
&=\sum_{n\in\mathbb Z}S_n^*z^{n}.\nonumber
\end{align}

It follows that 
for $n\geq 0$
$$
S_n.1=\delta_{n, 0}, \ \ S_{-n}^*.1=\delta_{n, 0}.
$$

There is a natural hermitian inner product on $V$ given by 
\begin{equation}
b_n^* = b_{-n}.
\end{equation}

The elements 
$b_{-\lambda}=b_{-\lambda_1}\cdots b_{-\lambda_l}
(\lambda \vdash n)$ span an orthogonal basis in $V$ and
\begin{equation}
(b_{-\lambda},b_{-\mu})=
\delta_{m,n}z_{\lambda},
\end{equation}
where $z_{\lambda}=\prod_{i\geq 1}^{}i^{m_i}m_i!$
for $\lambda=(1^{m_1}2^{m_2}\cdots)$.

The vertex space $V$ is isomorphic to the ring $\Lambda_{\mathbb Q}$
under the characteristic \linebreak
map:
$$
b_{-\la}=b_{-\la_1}\cdots b_{-\la_l}\longrightarrow p_{\la}.
$$
The characteristic map turns out to be also isometric (see \cite{J1}).

\begin{definition} The polynomial $s_{n}=S_{-n}.1$
is called the $n^{th}$ complete homogeneous symmetric polynomial in the $b_{-n}$.
\end{definition}

\begin{remark} We do not attempt to distinguish
the term "symmetric function" from that of "symmetric polynomial". It is apparant that the symmetry is
only with respect to the variables $x_1, x_2, \cdots$
and not with respect to the power sum variable.
\end{remark}

The generating function of $s_{n}$ is given by:

\begin{equation}
\sum_{n=0}^{\infty}s_{n}z^n = exp(\sum_{n=1}^{\infty}
\frac{b_{-n}}{ n}z^n).
\end{equation}

\begin{lemma} \cite{J2} The components of $S(z)$ and $S^*(z)$
satisfy the following commutation relations.
\begin{align*}
S_mS_n+S_{n-1}S_{m+1}&=0, \qquad S_m^*S_n^*+S_{n+1}^*S_{m-1}^*=0,\\
S_mS_n^*+S_{n-1}^*S_{m-1}&=\delta_{m,n}.
\end{align*}
\end{lemma}

Two $l$-tuple $\mu$ and $\la$ are related if there is a permutation $\sigma$ such that $\mu+\delta=\sigma(\la+\delta)$. If $\mu_i=\mu_{i+1}-1$, then $\mu+\delta=\sigma(\mu+\delta)$ for $\sigma=(i, i+1)$.
If there exists an old permutation $\sigma$ such that
$\mu+\delta=\sigma(\mu+\delta)$ then we say that $\mu$ is degenerate.
For any non-degenerate
$l$-tuple $\mu=(\mu_1, \cdots, \mu_l)$ in the general position
such that no two parts $\mu_i, \mu_j (i<j)$ 
satisfying $\mu_i=\mu_{j}-(j-i)$
 then there exists uniquely
a partition $\la$ such that $\mu+\delta=\sigma(\la+\delta)$. We remark  that the degeneracy condition corresponds to whether a weight belongs to
some wall of the Weyl chambers for $\mathfrak S_l$. In the latter context the symmetry is called the linkage symmetry. For simplicity we let $sgn(\mu)=0$ if the $l$-tuple $\mu$ is degenerate.
 
\begin{theorem} 
\cite{J1, J5} \label{T:theorem1}
Under the mapping $b_{-n}\rightarrow p_n$, the 
space $V$
is isometrically isomorphic to $\Lambda_{\mathbb Q}$.
The set $\{h_{-\lambda}=s_{\lambda_1}s_{\lambda_2}
\cdots s_{\lambda_l}:\lambda \vdash n, n\in \mathbb Z_+\}$ 
and
$\{S_{-\la_1}S_{-\la_2}\cdots S_{-\la_l}.1: \lambda \vdash n, n\in \mathbb Z_+\}$
both form ${\mathbb C}$-linear bases. Moreover the basis 
$\{S_{-\la_1}S_{-\la_2}\cdots S_{-\la_l}.1\}$ is 
orthonormal and expressed
 explicitly by:
\begin{equation}
S_{-\la}.1:=S_{-\la_1}S_{-\la_2}\cdots S_{-\la_l}.1=det(s_{\la_i-i+j})=s_{\la},
\end{equation}
where $s_{\lambda}$ is the Schur function in the $b_{-n}$ under the 
isomorphism
 between $\Lambda_{\mathbb Q}$ and ${V}$.
\end{theorem}
\begin{proof} To show the $\mathfrak S_n$-symmetry we consider the modified 
vertex operators associated with the root lattice $\mathbb Z\alpha$ with $(\alpha|\alpha)=1$.
 Let   $\tilde{V}=V\otimes \mathbb C[\mathbb Z\alpha]$, where $\mathbb 
C[\mathbb Z\alpha]$ is the group algebra generated by $e^{m\alpha}, 
m\in \mathbb Z$. Define
\begin{align}
\overline{S}(z)&=S(z)e^{\alpha}z^{\partial}=\sum_{n\in \mathbb Z+1/2}
\overline{S}_n z^{-n-1/2},\\
\overline{S}^*(z)&=S^*(z)e^{-\alpha}z^{-\partial}=\sum_{n\in \mathbb Z+1/2}
\overline{S}^*_n z^{n-1/2},
\end{align}
where the operators $e^{\alpha}$ and $z^{\partial}$ act
on $\mathbb C[\mathbb Z\alpha]$ as follows: 
\begin{equation*}
e^{\alpha}e^{m\alpha}=e^{(m+1)\alpha}, \qquad
z^{\partial}e^{m\alpha}=z^me^{m\alpha},
\end{equation*}
The components satisfy the Clifford algebra relations
\begin{align} \label{E:clifford1}
\{\overline{S}_m, \overline{S}_n\}&=\{\overline{S}_m^*, \overline{S}_n^*\}=0,\\ \label{E:clifford2}
\{\overline{S}_m, \overline{S}_n^*\}&=\delta_{m,n},
\end{align}
where $m, n\in \mathbb Z+1/2$.

For any degenerate $l$-tuple $\la=(\la_1, \la_2, \cdots, \la_l)$
satisfying $\la_i=\la_{j}-j+i$ for some $i$ and $j$ we have
\begin{align*}
S_{-\la_1}\cdots S_{-\la_l}.1&= 
\overline{S}_{-\la_1+1/2}\overline S_{-\la_2-1/2}
\cdots \overline{SF}_{-\la_l-l+1/2}e^{-l\alpha}\\
&=\overline S_{-\la+\delta-(m-1/2){\bf 1}}\\
&=0
\end{align*}
due to $\la+\delta=(i,j)(\la+\delta)$ and 
(\ref{E:clifford1}--\ref{E:clifford2}). Here $\bf 1$
denotes the integer vector $(1, 1, \cdots, 1)$, and similarly $m{\bf 1}=(m, m, \cdots, m)$.
In general for any non-degenerate
$l$-tuple we have
\begin{equation}
S_{-\mu_1}\cdots S_{-\mu_l}.1=sgn(\mu)S_{-\la_1}\cdots
S_{-\la_l}.1
\end{equation}
where $\la$ is the partition related to $\mu$: $\la=\sigma(\mu+\delta)-\delta$, and the sign $sign(\mu)=(-1)^{l(\sigma)}$.  See \cite{J1}  for details.
\end{proof}

The following result will be useful in our discussion.

\begin{proposition} \label{P:schur} 
Let $\delta=(n-1, n-2, \cdots, 1, 0)$. The
operator products \linebreak
$S(z_1)S(z_2)\cdots S(z_n)z^{\delta}$
and $S^*(z_1)S^*(z_2)\cdots S^*(z_n)z^{\delta}$
are skew-symmetric under the action of $\mathfrak S_n$. For any $w\in
\mathfrak S_n$ we have
\begin{align*}
S(z_{w(1)})S(z_{w(2)})\cdots S(z_{w(n)})z^{w(\delta)}
&=(-1)^{l(w)}S(z_1)S(z_2)\cdots S(z_n)z^{\delta}\\
S^*(z_{w(1)})S^*(z_{w(2)})\cdots S^*(z_{w(n)})z^{w(\delta)}
&=(-1)^{l(w)}S^*(z_1)S^*(z_2)\cdots S^*(z_n)z^{\delta}
\end{align*}
\end{proposition}

\begin{definition} For an $l$-tuple $\mu$ we define the 
Schur function $s_{\mu}$ to be the symmetric function
corresponding to $S_{-\mu_1}\cdots S_{-\mu_l}.1$ under the
characteristic map. We will simply write
\begin{equation*}
s_{\mu}=S_{-\mu}.1=S_{-\mu_1}\cdots S_{-\mu_l}.1\, .
\end{equation*}
\end{definition}

The following fact follows easily from Theorem 
\ref{T:theorem1} .
\begin{equation}
s_{\mu}=\left\{\begin{array}{ll}
0 & \text{$\mu$ is degenerate}\\
sgn(\sigma)s_{\la} &\text{$\la=\sigma(\mu+\delta)-\delta, \la\in \mathcal P$}\end{array}.\right.
\end{equation}

For convenience we denote the partition $\la$ associated to the
tuple $\mu$ by $\la=\pi(\mu)$.
 
We remark that the characteristic map can be defined over $\mathbb Z$ if we
do it on the basis of homogeneous polynomials or Schur polynomials.
 Correspondingly we have similar results for the dual vertex
operators $S^*(z)$. Here the
vector $S_{n}^*.1$ corresponds to the elementary 
symmetric function $(-1)^n e_{n}=(-1)^ns_{(1^n)}$.

For a partition $\la=(\la_1, \cdots, \la_l)$ we denote by $\la'=(\la_1', \cdots, \la_k')$ the dual partition, where
$\la_i'=Card\{j:\la_j\geq i\}$. We also denote by 
$(\la, \mu)$ the juxtaposition of two partitions or tuples.
Note that $(\la, \mu)$ is generally not a partition.

\begin{theorem} \label{T:theorem2}
 The set $\{S_{\la_1}^*\cdots S_{\la_l}^*.1 \mid\la \vdash n, n\in \mathbb Z_+\}$ also forms an
 orthonormal basis.
 \begin{equation}
S_{\la_1}^*S_{\la_2}^*\cdots S_{\la_l}^*.1=
 (-1)^{|{\la}|}s_{{\la}'}=(-1)^{|{\la}|}
det(s_{\la_i'-i+j}),
\end{equation}
where $|\la|=\la_1+\la_2+\cdots+\la_l.$
Moreover we have
\begin{align*}
S_{-\mu_1}\cdots S_{-\mu_l}S_{\nu_1}^*\cdots S_{\nu_k}^*.1
&=(-1)^{|\nu|}sgn(\la, \mu)s_{\pi(\mu, \nu')}\\
S_{\mu_1}^*\cdots S_{\mu_k}^*S_{-\nu_1}\cdots S_{-\nu_l}.1
&=(-1)^{|\mu|}sgn(\mu, \nu)s_{\pi(\mu, \nu')'}
\end{align*}
where $(\mu, \nu)$ is the juxtaposition of $\mu$ and $\nu$, and the associated
partition $\pi(\mu, \nu)$ is obtained by 
$\pi(\mu,\nu)=\sigma((\mu,\nu)+\delta)-\delta$ for some
$\sigma \in \mathfrak S_{l+k}$. 
\end{theorem}

As a consequence of Theorem \ref{T:theorem2} it follows that
$S_{-\la}.1=(-1)^{|\la|}S^*_{\la}.1
$.
In summary 
we have the following symmetry property on the parts. \begin{equation}
s_{\mu'}=\left\{\begin{array}{ll}
0 & \text{$\mu$ is degenerate}\\
sgn(\sigma)s_{\pi(\mu)'} &\text{$\pi(\mu)=\sigma(\mu+\delta)-\delta, 
\pi(\mu)\in \mathcal P$}\end{array}.\right.
\end{equation}
Note that we apply the action of 
$\mathfrak S_n$ first on the $n$-tuples and then follow 
by the duality.

\begin{example}  
\begin{align*}
S_{-1}S_2^*S_2^*.1&=S_{-1}S_{-2}S_{-2}.1=0, \\
S_{-1}S_2^*S_1^*S_1^*S_1^*.1&=-S_{-1}S_{-4}S_{-1}.1=S_{-3}
S_{-2}S_{-1}.1=s_{(3,2,1)}, 
\end{align*}
where $(1,4,1)+\delta=(3,5,1)\overset{-}{\sim} (5, 3,1 )=(3, 2, 1)+\delta$.
\end{example}

To close this section we derive the Littlewood-Richardson rule in our picture, which will be used later to realize
the action of $U_q(\widehat{sl}_2)$ on $\Lambda_{\mathbb Q[q, q^{-1}]}$. 

\begin{proposition} \label{P:littlewood}
Let $\lambda$ and $\mu$ be two partitions of lengths m and n
respectively. Then
$s_{\lambda}s_{\mu}=\sum{}{}C_{MN}s_{(\lambda+M,\mu-N)}$
where $C_{MN}$ is the number of integral matrices $(k_{ij})$ such that 
\begin{align*}
(k_{11}+k_{12}+\cdots+k_{1n},\cdots, k_{m1}+\cdots+k_{mn})&=M, \\
(k_{11}+k_{21}+\cdots+k_{m1},\cdots, k_{1n}+k_{2n}+\cdots+k_{mn})&=N.
\end{align*}
In particular the Schur functions form a $\mathbb Z$-basis
in $\Lambda$.
\end{proposition}
\begin{proof} It follows from definition that
\begin{equation*}
s_{\lambda}s_{\mu}=\int S(z_1)\cdots S(z_m).1S(w_1)\cdots S(w_n).1
z^{-\lambda}w^{-\mu}\frac{dz}{z}\frac{dw}{w}.
\end{equation*}

Observe that
\begin{align*}
&S(z_1)\cdots S(z_m).1 S(w_1)\cdots S(w_n).1\\
&=\prod_{i,j}^{}(1-\frac{w_j}{z_i})^{-1}S(z_1)\cdots S(z_m)S(w_1)\cdots
S(w_n).
\end{align*}

As an infinite series in $|w_j|<|z_i|$ we have 
\begin{align*}
\prod_{i,j}(1-\frac{w_j}{z_i})^{-1 }
&=\prod_{i,j}(1+\frac{w_j}{z_i}+(\frac{w_j}{z_i})^2+\cdots)\\
&=\prod_{i,j}(\sum_{k_{ij}}w_j^{k_{ij}}z_i^{-k_{ij}})\\
&=\prod_{k=(k_{ij})}w_1^{k_{\cdot 1}}\cdots w_n^{k_{\cdot
n}}z_1^{-k_{1\cdot}}\cdots z_m^{-k_{m\cdot}}, \end{align*}
with $k_{\cdot j}=k_{1j}+\cdots+k_{mj}$,
$k_{i\cdot}=k_{i1}+\cdots+k_{in},k_{ij}\geq0.$

Plugging the expansion into the integral and invoking Theorem (\ref{T:theorem1}).
 We prove the proposition.
\end{proof}

\begin{remark} The number $C_{MN}$ is equal to the index of $\mathfrak S_{\lambda}\pi \mathfrak S_{\mu}$
in $\mathfrak S_n$ \cite{JK}. We can also write:
$$
s_{\lambda}s_{\mu}=\prod_{i,j}(1-R_{ij})^{-1}s_{(\lambda,\mu)},
$$
\noindent where $R_{ij}$ is the raising operator defined on the parts of the Schur functions: $R_{ij}s_{(\cdots, \la_i, \cdots, \la_j, \cdots)}=s_{(\cdots, \la_i+1, \cdots, \la_j-1, \cdots)}$.
\end{remark}

\section{ Quantum affine algebra $U_q(\widehat{sl}_2)$.}
\label{S:quantum}

For $n\in\mathbb Z_+$ we define $[n]=\frac{q^n-q^{-n}}{q-q^{-1}}$.
The $q$-factorial $[n]!$ denotes $[n][n-1]\cdots [2][1]$
and then the $q$-Gaussian numbers are defined naturally by
$\left[\begin{matrix} n\\
  m\end{matrix}\right]=\frac{[n]!}{[m]![n-m]!}$
for $n\geq m\geq 0$. By convention $[0]=[1]=1$.
For an element $a$ in an algebra over $\mathbb Z[q, q^{-1}]$ we
use $a^{(n)}$ to denote the divided power $\frac{a^n}{[n]!}$.

The quantum affine algebra $U_q(\widehat{sl}_2)$ is the associative algebra
generated by Chevalley generators $e_i, f_i, K_i$ ($i=0, 1$) and $q^{d}$ subject to the following defining relations. 
\begin{align*}
  &K_iK_i^{-1}=K_i^{-1}K_i=1, \ \  q^dq^{-d}=q^{-d}q^d =1,\\
 &K_iK_j=K_jK_i,\ \  q^dK_i^{\pm 1}=K_i^{\pm 1}q^d,\\
 &K_ie_{j} K_i^{-1}=q^{ a_{ij}}e_{j},\ \ K_if_{j} K_i^{-1}=q^{-a_{ij}}f_{j},\\
&q^de_{i} q^{-d}=q^{ \delta_{i,0}}e_{i},\ \ 
q^df_{i} q^{-d}=q^{-\delta_{i,0}}f_{i},\\
&[e_{i}, f_{j}
]=\delta_{ij}\frac{K_i-K_i^{-1}}{q-q^{-1}},\\
 & \sum_{r=0}^{1-a_{ij}}(-1)^r
e_{i}^{(r)}e_{j}e_{i}^{(1-a_{ij}-r)}=0\
  \ \ \ \ \text{if $i\ne j$},\\
&\sum_{r=0}^{1-a_{ij}}(-1)^r f_{i}^{(r)}f_{j}f_{i}^{(1-a_{ij}-r)}=0\
  \ \ \ \ \text{if $i\ne j$}.
\end{align*}
where $(a_{ij})=\begin{pmatrix} 2 &-2\\-2 &2\end{pmatrix}$ is the 
extended Cartan matrix \cite{Ka}. 
The central element $K_0K_1=q^c$ acts as an integral power of $q$ on an integrable module of $U_q(\widehat{sl}_2)$. 

We will be interested in
level one ($c=1$) irreducible representations, which
had been realized by vertex operators in \cite{FJ}.
For our purpose we will recall the vertex representations
in a rescaled form.

Let $\{a_n|n\neq 0\}\bigcup\{c\}$ be the generators the Heisenberg algebra 
$U_q(\widehat{sl}_2 )$ with the defining relation.
\begin{equation} \label{E:heisenberg}
[a_m, a_n] =\delta_{m,-n}\frac{m}{1+q^{2|m|c}}.
\end{equation}

Following Frenkel-Jing \cite{FJ} the level one irreducible representation 
$V(\Lambda_i)$ is realized on the vertex representation space
$Sym({a_{-n}}'s)\bigotimes{\mathbb C}[{\mathbb Z}\alpha]e^{i\alpha/2}$ ($i=0, 1$)
at $c=1$. Here $Sym({a_{-n}}'s)$ denotes the symmetric algebra
generated by the Heisenberg generators
$a_{-n}$'s. The element $\vi$ (the highest weight vector) 
will be formally adjoined
to $\mathbb C[\mathbb Z\alpha]$.

We define two kinds of operators on the vector space
 ${\mathbb C}[{\mathbb Z}\alpha]e^{i\alpha/2}=<e^{m\alpha}e^{i\alpha/2}|\linebreak
m \in {\mathbb Z}>$:
\begin{align}
e^{n\alpha}.e^{m\alpha}e^{i\alpha/2} &= e^{(m+n)\alpha}e^{i\alpha/2}, n\in {\mathbb Z},\\
\partial.e^{m\alpha}e^{i\alpha/2}&=(2m+i)e^{m\alpha}\vi.
\end{align}
In particular 
${\mathbb C}[{\mathbb Z}\alpha]e^{i\alpha/2}$ is a
${\mathbb C}[{\mathbb Z}\alpha]$-module.

We define the vertex operators associated to $U_q(\widehat{sl}_2)$ by:
\begin{align}
X^+(z)&=exp(\sum_{n=1}^{\infty}\frac{(1+q^{2n})q^{-n}}{ n}a_{-n}z^n)
exp(-\sum_{n=1}^{\infty}\frac{(1+q^{2n})q^{-n}}{ n}a_nz^{-n})e^{\alpha}z^{\partial}\\
&=\sum_{n\in{\mathbb Z}}X_n^+z^{-n-1}\nonumber\\
X^-(z)&=exp(-\sum_{n=1}^{\infty}\frac{1+q^{2n}}{ n}a_{-n}z^n)
exp(\sum_{n=1}^{\infty}\frac{1+q^{2n}}{n}a_nz^{-n})e^{-\alpha}z^{-\partial}\\
&=\sum_{n\in{\mathbb Z}}X_n^-z^{-n-1}\nonumber
\end{align}

The normal order of vertex operator products is defined
by rearranging the exponential factors. For example,
\begin{align*}
&:X^+(z)X^+(w):\\
&=exp(\sum_{n=1}^{\infty}\frac{(1+q^{2n})q^{-n}}{ n}a_{-n}(z^n+w^n))
exp(-\sum_{n=1}^{\infty}\frac{(1+q^{2n})q^{-n}}{ n}a_n
(z^{-n}+w^{-n}))\\
&\qquad\qquad e^{2\alpha}z^{\partial}w^{\partial}.
\end{align*}

The components of the Drinfeld generators satisfy some
quadratic relations (Serre relations).

\begin{lemma} \cite{FJ} The components of $X^{\pm}(z)$ 
satisfy the following commutation relations.
\begin{align*}
X_m^{\pm}X_n^{\pm}-q^{\pm 2}X_n^{\pm}X_m^{\pm}&=
q^{\pm 2}X_{m-1}^{\pm}X_{n+1}^{\pm}-X_{n+1}^{\pm}X_{m-1}^{\pm},\\
X_m^{+}X_n^{-}-X_n^{-}X_m^{+}&=\frac 1{q-q^{-1}}
\left(\psi_{m+n}q^{(m-n)/2}-\phi_{m+n}q^{-(m-n)/2}\right),
\end{align*}
where the polynomials $\psi_n$ and $\phi_{-n}$ are defined by
\begin{align}
\Psi(z)&=\sum_{n\geq 0}\psi_n z^{-n}=
exp(\sum_{n\in \mathbb N}\frac {(q^{2n}-q^{-2n})q^{n/2}}n a_n z^{-n})q^{\partial}\\
\Phi(z)&=\sum_{n\geq 0}\phi_{-n} z^{n}=
exp(\sum_{n\in \mathbb N}\frac {(q^{-2n}-q^{2n})q^{n/2}} n a_{-n} z^{-n})q^{-\partial}
\end{align}
\end{lemma}

With the action of $X^{\pm}_m$ the Chevalley generators are expressed by
\begin{alignat*}{3}
e_1 &\rightarrow X_0^+, &f_1&\rightarrow X_0^-, &K_1&\rightarrow q^\partial\\
e_0 &\rightarrow X_1^-q^{-\partial},& f_0&\rightarrow q^\partial X_{-1}^+,
& K_0&\rightarrow q^{1-\partial}
\end{alignat*}
The vertex space is endowed with the standard inner product via
\begin{align*}
a_n^* &= a_{-n}, \quad
(e^{\alpha})^* = e^{-\alpha},\\
(z^{\partial})^* &= z^{-\partial}.
\end{align*}

It follows from the commutation relations (\ref{E:heisenberg})
that 
\begin{equation}\label{E:form2}
(a_{-\lambda}e^{m\alpha}\vi, a_{-\mu}e^{n\alpha}\vi)=
\delta_{mn}\delta_{\lambda\mu}z_{\lambda}\prod_{j\geq 1}^{}
\frac{ 1}{1+q^{2\lambda_j}}
\end{equation}
where $z_{\lambda}$ as in section \ref{S:schur}.

By Section \ref{S:schur} there are two special bases in $V(\Lambda_i)$:
the power sum basis \linebreak
$\{a_{-\lambda}e^{m\alpha}\vi\}$ and the Schur basis 
$\{s_{\lambda}e^{m\alpha}\vi\}$.
However the Schur basis is no longer orthogonal with respect to the inner
product (\ref{E:form2}).

Let $b_{-n}=a_{-n}, b_n=(1+q^{2n})a_n$, $n\in {\mathbb N}$,
then
$\{b_{-n}\}$ generate a standard Heisenberg algebra
as in Section \ref{S:schur}. In terms of the new Heisenberg generators
we have
\begin{align*}
S(z)&=exp(\sum_{n=1}^{\infty}\frac{ 1}{n}a_{-n}z^n)
exp(-\sum_{n=1}^{\infty}\frac{1+q^{2n}}{ n}a_nz^{-n})
= \sum_{n\in{\mathbb Z}}^{}S_nz^{-n},\\
S^*(z)&=exp(-\sum_{n=1}^{\infty}\frac{ 1}{n}a_{-n}z^n)
exp(\sum_{n=1}^{\infty}\frac{1+q^{2n}}{ n}a_nz^{-n})
= \sum_{n\in{\mathbb Z}}^{}S_n^*z^{n}, 
\end{align*}
which generate the Schur function basis. 
For a partition $\la$ and $m\in\mathbb Z$ we define the {\it Schur symmetric polynomial} in $V_i$ (cf. Theorem \ref{T:theorem1}):
\begin{equation}
s_{\lambda}e^{m\alpha}\vi:=S_{-\lambda}e^{m\alpha}\vi=S_{-\lambda_1}\cdots S_{-\lambda_l}e^{m\alpha}\vi.
\end{equation}
The element $s_{\la}$ is a 
polynomial over $\mathbb Q$ in terms of the power sum $a_{\mu}$, where
$|\mu|=|\la|$. Note that $S(z)$ acts trivially on the 
lattice vector $e^{m\alpha}\vi$.

We need to recall some further terminology about partitions.
Let $\la$ and $\mu$ be two partitions we write $\la \supset
\mu$ if the Young diagram of $\la$ contains that of $\mu$. The set difference $\la-\mu$ is called a skew diagram. The conjugate of a skew diagram $\theta=\la-\mu$ is
$\theta'=\la'-\mu'$ and we define
\begin{equation}
|\theta|=\sum \theta_j=|\la|-|\mu|.
\end{equation}
A skew diagram $\theta$ is a horizontal $n$-strip (resp.
a vertical $n$-strip) if $|\theta|=n$ and $\theta'_j\leq 1$
(resp. $\theta_j\leq 1$) for each $j$. Thus a horizontal (resp. vertical) strip has at most one column (resp. rows) in its diagram. We will also denote
by $(\la, \mu)$ the juxtaposition of $\la$ and $\mu$. 
As before $\pi(\la, \mu)$ is the partition
associated to $(\la, \mu)$ by performing the linkage
symmetry.  

For a partition $\mu$ and an integer $m$ we let $\mathcal V_n=\mathcal
V_n(m, \mu)$ be the set of the partitions of $\la$
such that the skew diagram $\la-(m, \mu')'$ is a vertical $n$-strip.
We also let $\mathcal H_n=\mathcal H_n(m, \mu)$ be the set of partitions $\la$
such that the skew diagram $\la-(m, \mu)$ is a horizontal $n$-strip. Note that $\mathcal V_n$ may be described as the set of partitions $\la$ such that the skew diagram $\la-(1^m, \mu)$ is a vertical $n$-strip. The following is called the Pieri rule \cite{M}:
\begin{align}
s_nS_{-m}S_{-\mu}.1 &=\sum_{\la\in \mathcal H_i(m, \mu)}sgn(m,\mu)s_{\la},\\
s_{1^n}S_{-m}S_{-\mu}.1 &=\sum_{\la\in \mathcal V_n(m, \mu)}(-1)^n 
sgn{(m, \mu')}'
s_{\la}.
\end{align}

We will also frequently use
$\la-\mu$ to denote the difference of two
integral vectors in $\mathbb Z^n$. 

\begin{theorem}\label{theorem3}
The quantum affine algebra $U_q(\widehat{sl}_2 )$ is realized
on the Fock space $\Lambda\bigotimes{\mathbb C}[{\mathbb Z}\alpha]\vi$ 
of symmetric functions by the following action:
\begin{align*}
X_n^+s_{\mu}e^{m\alpha}\vi
&=\left(\sum_{j=0}^{l(\mu)-2m-n-1-i}q^{-2m-n-1-i-2j}\cdot\right.
\\&\left.  sgn(-2m-n-1-i-j,\mu)
\sum_{\lambda\in \mathcal H_j}
s_{\lambda}\right)e^{(m+1)\alpha}\vi.
\end{align*}
where $\mathcal H_j=\mathcal H_j(-2m-n-1-i-j, \mu)$, 
the sign refers to $sgn(-2m-n-1-i-j, \mu)=(-1)^{l(\sigma)}$
such that 
$(-2m-n-1-i-j, \mu)+\delta=\sigma(\la+\delta)$
and $\la$ is the partition of length at most $l(\mu)+1$.
Also we have
\begin{align*}
X_n^-s_{\mu}e^{m\alpha}\vi
&=(-1)^{n+1+i}\left(\sum_{j=0}^{\mu_1+2m-n-1+i}q^{2j}\cdot\right.\\
&\left.sgn(2m-n-1-j+i, \mu')\sum_{\lambda\in\mathcal V_j }s_{\lambda}\right)
e^{(m-1)\alpha}\vi
\end{align*}
where $\mathcal V_j=\mathcal V_j(2m-n-1-j+i, \mu')$,
the sign refers to $sgn(2m-n-1-i-j, \mu)=(-1)^{l(\sigma)}$
such that 
$(2m-n-1-j+i, \mu)+\delta=\sigma(\la+\delta)$
and $\la$ is the partition of length at most $l(\mu)+2m-n-1-j+i$.
\end{theorem}

\begin{proof} This follows from our vertex operator calculus of symmetric functions
\cite{J1, J2}.
\begin{align*}
&X^+(z)S(w_1)\cdots S(w_l)e^{m\alpha}\vi\\
&=:X^+(z)S(w_1)\cdots S(w_l):z^{2m+i}e^{(m+1)\alpha}\vi
\prod_{j}^{}(1-q^{-1}\frac{w_j}{z})\prod_{j<k}^{}(1-\frac{w_k}{w_j})\\
&=exp(\sum_{n=1}^{\infty}\frac{q^{-2n}}{ n}a_{-n}w_0^n)
S(w_0)\cdots S(w_l)e^{(m+1)\alpha}\vi w_0^{2m+i}q^{-2m-i},
\end{align*}
where $w_0=qz$. Taking the coefficient of $z^{-n-1}w^{-\mu}$ and using Littlewood-Richard\-son rule (\ref{P:littlewood}) we
obtain the result. The case of $X_n^-s_{\mu}e^{m\alpha}$
is proved similarly.
\begin{align*}
&X^-(z)S(w_1)\cdots S(w_l)e^{m\alpha}\vi\\
&=:X^-(z)S(w_1)\cdots S(w_l):z^{-2m-i}e^{(m-1)\alpha}\vi
\prod_{j<k}(1-\frac {w_k}{w_j})
\prod_{j}(1-\frac {w_j}{z})^{-1}\\
&=S^*(q^2z).1S^*(z)S(w_1)\cdots S(w_l)z^{-2m-i}e^{(m-1)\alpha}\vi.
\end{align*}
Taking the coefficient of $z^{-n-1}w^{\mu}$ we obtain the formula.
\end{proof}

We can reformulate the result in terms of the standard
inner product in Section \ref{S:schur}. Let 
$u_j=(0, \cdots, 0, 1, 0, \cdots, 0)$ be the $j$th unit vector
in $\mathbb Z^n$.
Let ${\bf 1}_{(l_1, \cdots, l_j)}$ be the sum of the unit vectors
$u_{l_1}, \cdots, u_{l_j}$.

\begin{proposition} For $n\in \mathbb Z$ and a partition
$\mu$ we have
\begin{align*}
&X_n^-s_{\mu}e^{m\alpha}\vi\\
&=\sum_{\la}s_{\la}e^{(m-1)\alpha}\vi\sum_{j=0}^{l(\la)}
(-q^2)^j
\sum_{l_1<\cdots <l_j}(S_{2m+i-n-1-j}^*S_{-\mu},
S_{-(\la-{\bf 1}_{(l_1, \cdots, l_j)})}),
\end{align*}
where $\la$ runs through partitions of weight
$|\mu|+2m-n-1+i$ such that $\la-{\bf 1}_{(l_1, \cdots, l_j)}$
is the juxtaposition of $(1^{2m-n-j-1+i})$ and $\mu$.
\end{proposition}

Later in Theorem \ref{T:dividedpower} we will give another
proof in terms of the dual vertex operator $S^*(z)$.

\begin{example} Using Theorem \ref{theorem3} it is easy to compute
the following.

\begin{align*}
&X_n^{\pm}e^{r\alpha}\vi=0, \qquad\mbox{if $n>\mp 2r-1\mp i$},\\
&X^+_{-2r+1-i}\cdots X^+_{-3-i}X^+_{-1-i}\vi=e^{r\alpha}\vi, \qquad r\geq 1, \\
&X^-_{2r-3+i}\cdots X^+_{1+i}X^+_{-1+i}\vi=e^{-r\alpha}\vi,
\qquad r\geq 1.
\end{align*}
\begin{align*}
X_{-1}^+s_{(2,1)}e^{-\alpha}&=
q^2s_{(2,2,1)}-
q^{-2}(s_5+s_{(4,1)}+s_{(3, 2)})+
q^{-6}(s_5+s_{(4,1)}),\\
X_0^-s_{1}e^{\alpha}&=-s_{(2)}+q^4s_{(1^2)}.
\end{align*}
\end{example}

We can generalize the action to the divided powers of
$X_n^{\pm (r)}$.

\begin{lemma} \cite{BFJ} \label{L:qsymmetry}
For $r\in \mathbb N$ we have
\begin{equation*}
\prod_{1 \le i < j \le k} (z_i - qz_j) =
\sum_{w \in \mathfrak S_k} (-q)^{\ell(w)} z^{w(\delta)} +
 \sum a_{\gamma_1, \dots, \gamma_k}
 z_1^{\gamma_1}  z_2^{\gamma_2} \dots z_k^{\gamma_k},
\end{equation*}
where $\delta = (k-1, k-2,\dots, 0)$
and the second sum consists of certain monomials such that some $\gamma_i =
\gamma_j$, $i \neq j$ and 
$a_{{\gamma}} \in \mathbb Z[q],  \ 
a_{{\gamma}}(1) = 0$.
\end{lemma}

\begin{theorem}\label{T:dividedpower}
The lattices $\mathcal L_i=\sum_{\la, m}\oplus \mathbb Z[q, q^{-1}]s_{\la}e^{m\alpha}\vi$ are invariant under the action of the divided powers ${X_{n}^\pm}^{(r)}$. More precisely we have
\begin{align*}
&X_n^{+(r)}s_{\mu}e^{m\alpha}\vi\\
&=q^{-3\binom{r}{2}-r(n+1+2m+i)}\sum_{l(\la)\leq r}
q^{-2|\la|}s_{\la}\cdot
s_{(-\la-2\delta-(n+1+2m+i){\bf 1}, \mu)}e^{(m+r)\alpha}\vi\\
&X_n^{-(r)}s_{\mu'}e^{m\alpha}\vi\\
&=(-1)^{r(n+1+i)}q^{\binom{r}{2}}
\sum_{l(\la)\leq r}q^{2|\la|}
s_{\la'}\cdot
s_{(-\la-2\delta-(n+1-2m-i){\bf 1}, \mu)'}
e^{(m-r)\alpha}\vi
\end{align*}
where the summations run through all partitions $\la$
of length $\leq r$,   
$\delta =(r-1, r-2, \cdots, 0)$,
and ${\bf 1}=(1, \cdots, 1)\in \mathbb Z^r$. 
\end{theorem}
\begin{proof} Let $z=(z_1, \cdots, z_r)$, $w=(w_1, \cdots,  w_l)$, and ${\bf 1}=(1, \cdots, 1)\in \mathbb Z^r$ 
in the following computation.
\begin{align*}
    &{x_n^+}^r s_{\mu}e^{m \alpha}\vi\\
&= \oint X^+(z_1)\cdots X^+(z_r)S(w_1)\cdots S(w_l) z^{(n+1)\bf 1}w^{-\mu}  \ \frac{dzdw}{zw} \  e^{m \alpha}\vi  \\
& =  \oint \exp\bigl(\sum_{n=1}^\infty \frac{(q^n +
  q^{-n})}{n} a_{-n} (z_1^n + \dots + z_r^n)\bigr) 
:S(w_1)\cdots S(w_l):\\ & \hskip
 .2 in \times \prod_{i < j}
(z_i - z_j)(z_i - q^{-2}z_j)(1-\frac{w_j}{w_i})
\prod_{i, j} (1-q^{-1}\frac {w_j}{z_i})\\
&\qquad\qquad \times z^{(2m + n+1+i){\bf 1}} e^{(m + r)\alpha}\vi  \frac{dz}z\frac{dw}{w}.
\end{align*}
 Note that the integrand divided by $\prod_{j<k} (z_j - q^{-2} z_k)$ is an
anti-symmetric function in $z_1, \dots, z_r.$ It follows from Lemma \ref{L:qsymmetry}
that the terms $z_1^{\gamma_1} z_2^{\gamma_2} \dots
z_r^{\gamma_r}$ (for which some $\gamma_k = \gamma_j$) make no contribution
to the integral.  Therefore
\begin{align*} 
x_n^{+(r)} & s_{\mu}e^{m \alpha}\vi   = \frac 1{[r]!}\sum_{w \in \mathfrak S_r} 
\oint  \exp \bigl( \sum_{n=1}^\infty
  \frac{(q^n + q^{-n})}{n}
  a_{-n} (z_1^n + \dots + z_r^n)\bigr)\\
&\times :S(w_1)\cdots S(w_l):\prod_{i < j} (z_i - z_j)
\prod_{i < j}(1-\frac{w_j}{w_i})
\prod_{i, j} (1-q^{-1}\frac {w_j}{z_i})\\
 & (-q)^{-\ell(w)} z^{w(\delta) + (2m + n + 1+i)
    {{\bf 1}}}w^{-\mu} e^{(m+r) \alpha}\vi
 \frac{dzdw}{zw} \\ 
& = q^{-\binom r2}
 \oint \exp\bigl( \sum_{n = 1}^\infty
  \frac{q^{-2n}}{n} a_{-n} (z_1^n + \dots + z_r^n)
  \bigr) S(qz_1)\cdots S(qz_r)\\
&S(w_1)\cdots S(w_l)
z^{2\delta + (2m + n + 1+i){{\bf 1}}}w^{-\mu} e^{(m + r) \alpha}\vi
\frac{dzdw}{zw} ,
\end{align*}
where $\delta = (r-1, r-2, \dots, 0)$ and we have used 
$$\sum_{w \in\mathfrak S_r} q^{-2\ell(w)} = q^{-\binom{r}{2}} [r]!.$$  
From the orthogonality of Schur functions \cite{M} it follows that
\begin{equation*}
\exp\bigl( \sum_{n = 1}^\infty
  \frac {a_{-n}}{n} (z_1^n + \dots + z_r^n)\bigr)
=\sum_{l(\la)\leq r} s_{\la}(a_{-k})s_{\la}(z_i)
\end{equation*}
where $s_{\la}(a_{-k})$ is the Schur function
in terms of the power sum $a_{-\mu}$ and $s_{\la}(z_i)$ is the
Schur polynomial in the variables $z_1, \cdots, z_r$. 
Replacing $z_j$ by
$q^{-1}z_j$ we get
\begin{align*}
&{x_n^+}^{(r)}s_{\mu} e^{m \alpha}\vi \\ 
&=q^{- 3\binom{r}{2}-r(n+1+2m+i)}
\sum_{l(\la)\leq r}q^{-2|\la|}s_{\la}\oint s_{\la}(z)z^{2\delta+(n+1+2m+i){\bf 1}} \\
& \times S(z_1)\cdots S(z_r)S(w_1)\cdots S(w_l)e^{(m+r)\alpha}\vi
w^{-\mu}\frac {dz}{z}\frac {dw}{w}\\
&=q^{-3\binom{r}{2}-r(n+1+2m+i)}\sum_{l(\la)\leq r}s_{\la}\oint \sum_{w\in\mathfrak S_r}\frac{(-1)^{l(w)}z^{w(\la+\delta)}}{\prod_{j<k}(z_j-z_k)}
z^{\delta+(n+1+2m+i){\bf 1}}\\
& \times z^{\delta}S(z_1)\cdots S(z_r)S(w_1)\cdots S(w_l)e^{(m+r)\alpha}\vi
w^{-\mu}\frac {dz}{z}\frac {dw}{w}.
\end{align*}
Applying the symmetry of the Schur vertex operators
in Proposition \ref{P:schur} we see that the above
expression becomes
\begin{align*}
&q^{- 3\binom{r}{2}-r(n+1+2m+i)}\sum_{l(\la)\leq r}q^{-2|\la|}s_{\la}
\sum_{w\in\mathfrak S_r}\oint\frac{(-1)^{l(w)}z^{\la+\delta}}{\prod_{j<k}(z_j-z_k)}
z^{w(\delta)+(n+1+2m+i){\bf 1}}\\
& \times z^{\delta}S(z_1)\cdots S(z_r)S(w_1)\cdots S(w_l)e^{(m+r)\alpha}\vi
w^{-\mu}\frac {dz}{z}\frac {dw}{w}\\
&=q^{- 3\binom{r}{2}-r(n+1+2m+i)}\sum_{l(\la)\leq r}q^{-2|\la|}s_{\la}\oint
z^{\la+2\delta+(n+1+2m+i){\bf 1}}\\
& \times S(z_1)\cdots S(z_r)S(w_1)\cdots S(w_l)e^{(m+r)\alpha}\vi
w^{-\mu}\frac {dz}{z}\frac {dw}{w},
\end{align*}
where we have used the Weyl denominator formula (see $\la=0$ in (\ref{E:weyl}))
and the integral is taken along contours in $z_i, w_i$ around the origin. The formula for $X_n^{+(r)}$
is then obtained by using Theorem \ref{T:theorem1}.

 The case of $X_n^{-(r)}$ is proved similarly with the help of the dual vertex operator $S^*(z)$. 
\begin{align*}
    &{X_n^-}^{(r)} s_{\mu'}e^{m \alpha}\vi\\
&= \frac {(-1)^{|\mu|}}{[r]!}\oint X^-(z_1)\cdots X^-(z_r)S^*(w_1)\cdots S^*(w_l) z^{n+1} 
w^{-\mu}e^{(m-r) \alpha}\vi\frac{dzdw}{zw}  \\
&=\frac {(-1)^{|\mu|}}{[r]!}\oint :S^*(q^2z_1)\cdots S^*(q^2z_r):S^*(z_1)\cdots S^*(z_r)S^*(w_1)\cdots S^*(w_l)\\
&\qquad \times \prod_{i< j}(z_i-q^2z_j)z^{\delta+(n+1-2m-i){\bf 1}}w^{-\mu} \  e^{(m-r) \alpha}\vi\frac{dzdw}{zw}\\
&=q^{\binom r2}\oint
:S^*(q^2z_1)\cdots S^*(q^2z_r):S^*(z_1)\cdots S^*(z_r)S^*(w_1)\cdots S^*(w_l)\\
&\times z^{2\delta+(n+1-2m-i){\bf 1}}w^{-\mu} \  e^{(m-r) \alpha}\vi\frac{dzdw}{zw}(-1)^{|\nu|}
\end{align*}
where we have used the skew-symmetry of
the integrand and Lemma \ref{L:qsymmetry}. Then the
formula for $X_n^{-(r)}$ is obtained from the following identity and the Weyl denominator formula.
\begin{equation*}
\exp\bigl( -\sum_{n = 1}^\infty
  \frac {a_{-n}}{n} (z_1^n + \dots + z_r^n)\bigr)
=\sum_{l(\la)\leq r} s_{\la'}(a_{-k})s_{\la}(z_i)
\end{equation*}

\end{proof}


\section{Combinatorial realization of $U_q(\widehat{sl}_2)$}\label{S:comb}

Let $\mathcal A=\mathbb Z[q, q^{-1}]$. We denote by $U_q(\widehat{sl}_2)_{\mathcal A}$ the
$\mathcal A$-algebra generated by
$e_i^{(n)}$, $f_i^{(n)}$, $K_i^{\pm 1}$. From section \ref{S:quantum}
it follows that

\begin{align}
&e_1^{(r)}s_{\mu}e^{m\alpha}\vi
=q^{-3\binom{r}{2}-r(2m+i+1)}\nonumber\\ \label{E:comb1}
&\times\left(\sum_{l(\la)\leq r}
q^{-2|\la|}s_{\la}
s_{(-\la-2\delta-(2m+1+i){\bf 1}, \mu)}\right)e^{(m+r)\alpha}\vi\\
&f_1^{(r)}.s_{\mu'}e^{m\alpha}\vi
=(-1)^{r(1+i)}q^{\binom{r}{2}}\nonumber\\ \label{E:comb2}
&\times\left(\sum_{l(\la)\leq r}
q^{2|\la|}s_{\la'}
s_{(-\la-2\delta+(2m+i-1){\bf 1}, \mu)'}\right)e^{(m-r)\alpha}\vi
\end{align}
\begin{align}
&f_0^{(r)}s_{\mu}e^{m\alpha}\vi
=q^{r(5-r)/2}\nonumber\\ \label{E:comb3}
&\times\left(\sum_{l(\la)\leq r}
q^{-2|\la|}s_{\la}
s_{(-\la-2\delta-(2m+i){\bf 1}, \mu)}\right)e^{(m+r)\alpha}\vi\\
&e_0^{(r)}.s_{\mu'}e^{m\alpha}\vi
=(-1)^{ri}q^{-\binom{r}{2}-r(2m+i)}\nonumber\\ \label{E:comb4}
&\times\left(\sum_{l(\la)\leq r}
q^{2|\la|}s_{\la'}
s_{(-\la-2\delta+(2m+i-2){\bf 1}, \mu)'}\right)e^{(m-r)\alpha}\vi
\end{align}

 As a consequence of these formulas and the Littlewood-Richardson rule (\ref{P:littlewood}) we get
the following theorem.

\begin{theorem}  The $\mathcal A$-lattice generated by the Schur functions $\{s_{\mu}e^{m\alpha}\}$, $m\in \mathbb Z, \mu\in\mathcal P$, 
is invariant under the integral form 
$U_q(\widehat{sl}_2)_{\mathcal A}$.
\end{theorem}

\begin{proposition} For $m\geq 0$ we have
\begin{align*}
f_1^{(2m)}e^{m\alpha}&=(-1)^mq^{m(2m-1)}e^{-m\alpha},\\
f_0^{(2m+1)}e^{-m\alpha}&=(-1)^mq^{-(2m+1)(m-2)}e^{-(m+1)\alpha},\\
f_0^{(2m)}e^{-m\alpha}e^{\alpha/2}&=(-1)^mq^{-m(2m-5)}e^{m\alpha},\\
f_1^{(2m+1)}e^{m\alpha}e^{\alpha/2}&=(-1)^mq^{m(2m+1)}
e^{-(m+1)\alpha}.
\end{align*}
\end{proposition}
\begin{proof} The four formulas are proved similarly.
Take $f_1^{(2m)}e^{m\alpha}$. Observe that
$-2\delta+(2m-1){\bf 1}=(-2m+1, -2m+3, \cdots, 2m-3, 2m-1)$
is of weight zero, thus only $\la=0$ contributes
to the summation in $f_1^{(2m)}e^{m\alpha}$ (see (\ref{E:comb2})). Since the longest element in $\mathfrak S_{2m}$ has inversion number $m(m-1)$, the sign of
$s_{-2\delta+(2m-1){\bf 1}}$ is $(-1)^{m(m-1)}=(-1)^m$.
\end{proof}

\begin{cor} For $m\geq 0$ we have
\begin{align*}
f_1^{(2m)}f_0^{(2m-1)}\cdots f_1^{(2)}f_0.1&=(-1)^mq^{3m^2}
e^{-m\alpha},   \\
f_0^{(2m+1)}f_1^{(2m)}\cdots f_1^{(2)}f_0.1&=q^{(m+1)(m+2)}
e^{(m+1)\alpha},\\
f_0^{(2m)}f_1^{(2m-1)}\cdots f_0^{(2)}f_1e^{\alpha/2}&=(-1)^mq^{m(m+2)}
e^{m\alpha}e^{\alpha/2},   \\
f_1^{(2m+1)}f_0^{(2m)}\cdots f_0^{(2)}f_1e^{\alpha/2}&=q^{3m(m+1)}
e^{-(m+1)\alpha}e^{\alpha/2},\\
\end{align*}
\end{cor}

\begin{example} In the following we abbreviate
$f_{i_1}^{(n_1)}\cdots f_{i_r}^{(n_r)}.1=f_{i_1}^{(n_1)}\cdots f_{i_r}^{(n_r)}$
in the basic representation $V(\Lambda_0)$.
\begin{align*}
f_0&=q^2e^{\alpha}\\
f_1f_0&=-q^2(1+q^2)s_1\\
f_1^{(2)}f_0&=-q^3e^{-\alpha}\\
f_0f_1f_0&=q^2(q^2+1)s_1e^{\alpha}\\
f_1f_0f_1f_0&=-(q^4+q^2)(s_2-q^4s_{1^2})\\
f_0f_1^{(2)}f_0&=-q^3(s_{1^2}+[3]s_2)\\
f_0f_1f_0f_1f_0&=q^2(1+q^2)^2(s_2+s_{1^2})e^{\alpha}\\
f_0^{(2)}f_1^{(2)}f_0&=q^4(s_2+[3]s_{1^2})e^{\alpha}\\
f_0^{(3)}f_1^{(2)}f_0&=q^6e^{2\alpha}\\
f_1^{(2)}f_0f_1f_0&=q^5(1+q^2)s_1e^{-\alpha}
\end{align*}

\end{example}

\begin{example} As in the last example we use 
$f_{i_1}^{(n_1)}\cdots f_{i_r}^{(n_r)}$ to denote \linebreak
$f_{i_1}^{(n_1)}\cdots f_{i_r}^{(n_r)}e^{\alpha/2}$ 
in the basic representation $V(\Lambda_1)$.
\begin{align*}
f_1&=e^{-\alpha}\\
f_0f_1&=q^{-2}(1+q^{-2})s_1\\
f_0^{(2)}f_1&=-q^3e^{\alpha}\\
f_1f_0f_1&=-(q^{-2}+1)s_1e^{-\alpha}\\
f_0f_1f_0f_1&=(1+q^2)(s_2-s_{1^2})\\
f_1f_0^{(2)}f_1&=-q^5([3]s_{1^2}+s_2)\\
f_1f_0f_1f_0f_1&=(1+q^2)^2(s_2+q^2s_{1^2})e^{-\alpha}\\
f_1^{(2)}f_0^{(2)}f_1&=q^5([3]s_2+s_{1^2})e^{-\alpha}\\
f_1^{(3)}f_0^{(2)}f_1&=q^6e^{-2\alpha}\\
f_0^{(2)}f_1f_0f_1&=-[2]s_1e^{\alpha}
\end{align*}
\end{example}

It would be interesting to see the relation between our formulas
and the fermionic picture \cite{LLT}.


\begin{thebibliography}{DJKM}

\bibitem[BCP]{BCP} J.~Beck, V.~Chari and A.~Pressley,
An algebraic characterization of the affine canonical basis, Duke Math. J., to
appear, math.QA/9808060.

\bibitem[BFJ]{BFJ} J.~Beck, I.B.~Frenkel and N.~Jing, Canonical basis and Macdonald polynomials,
Adv. in Math. {\bf 140} (1998), 95-127.

\bibitem[CP]{CP} V.~Chari and A.~Pressley,
 Finite dimensional
  representations of quantum affine algebras,
Representation Theory {\bf 1} (1997), 280--328.

\bibitem[DJKM]{DJKM} Date, Jimbo, Kashiwara and Miwa,
Transformation groups for soliton equations.
Nonlinear integrable systems---classical theory and quantum theory 
(Kyoto, 1981), pp. 39--119, World Sci. Publishing, Singapore, 1983.

\bibitem[F1]{F1} I.B.~Frenkel, Two constructions of affine Lie algebra 
representations and boson-fermion correspondence in quantum
field theory, J. Funct. Anal. {\bf 44} (1981), 259--327.

\bibitem[F2]{F2} I.B.~Frenkel, Lectures at Yale University, 1986.

\bibitem[FJ]{FJ} I.B.~Frenkel and N.~Jing, Vertex representations of quantum affine 
algebras, Proc. Natl. Acad. Sci. USA {\bf 85} (1988), 9373-9377.

\bibitem[G]{G} H.~Garland, The arithmetic theory of loop groups, J. Algebra {\bf 53} (1978), 480-551. 

\bibitem[H]{H} T.~Hayashi, Q-analogues of Clifford and Weyl algebras-spinor
and oscillator representations of quantum enveloping algebras,
Commun. Math. Phys. {\bf 127} (1990), 129-144.

\bibitem[JK]{JK} G.~James and A.~Kerber, {\em The representation theory of the symmetric group}, Encyclopedia of Math. and its
Appl. {\bf 16}
, Addison-Wesley, Reading, MA, 1981.

\bibitem[J1]{J1} N.~Jing, Vertex operators, symmetric functions and the spin groups $\Gamma_n$,
J. Algebra {\bf 138} (1991), 340-398.

\bibitem[J2]{J2} N.~Jing, Vertex operators and Hall-Littlewood functions, Adv. in Math.
{\bf 87} (1991), 226-248.

\bibitem[J3]{J3} N.~Jing, Vertex operators and generalized symmetric functions, in: Proc. of
Conf. on Quantum Topology (KSU, March 1993), ed. D.~Yetter, World Scientific, Singapore, 1994,
pp. 111-126.

\bibitem[J4]{J4} N.~Jing, q-Hypergeometric series and Macdonald functions, J. Alg. Comb.
{\bf 3} (1994) 291-305.

\bibitem[J5]{J5} N.~Jing, Boson-fermion correspondence for Hall-Littlewood polynomials,
J. Math. Phys. {\bf 36} (1995), 7073-7080.


\bibitem[J6]{J6} N.~Jing, Vertex representations of the quantum 
Kac-Moody algebras, Lett. Math. Phys. 
{\bf 44} (1998), no. 4, 261--271.

\bibitem[K]{K} M.~Kashiwara, Global crystal bases of quantum groups,
Duke Math. J. {\bf 73} (1993), 383--413.

\bibitem[Ka]{Ka} V.G.~Kac, {\em Infinite dimensional Lie algebras},
3rd. ed., Cambridge Univ. Press, Cambridge, 1990.

\bibitem[LLT]{LLT} A.~Lascoux, B.~Leclerc, J.Y.~Thibon, 
Ribbon tableaux, Hall-Littlewood functions, quantum affine algebras and 
unipotent varieties, S\'em. Lothar. Combin. {\bf 34} (1995), 23 pp. 

\bibitem[LP]{LP} J.~Lepowsky and P.~Primc, Structure of the standard modules for affine Lie algebra $A_1^{(1)}$,
Contemp. Math. {\bf 46}, Amer. Math. Soc., Providence, RI, 
1985.

\bibitem[L]{L} G.~Lusztig, {\em Introduction to quantum groups}, Progress in Mathematics {\bf 110}, Birkh\"auser, 
Boston, 1993.


\bibitem[M]{M} I.G.~Macdonald, {\em Symmetric functions and
Hall polynomials}, 2nd ed., Oxford University Press, New York, 1995.

\bibitem[MM]{MM} K.C.~Misra and T.~Miwa, Crystal bases for the basic
representations of $U_q(\hat{sl}(n))$, Commun. Math. Phys. {\bf 134} (1990), 79-88.

\bibitem[Z]{Z} A.~Zelevinsky, {\em Representations of finite classical 
groups, 
 A Hopf algebra approach.} Lecture Notes in
Mathematics, 869, Springer-Verlag, Berlin-New York, 1981. 


\end{thebibliography}
\end{document}